\documentclass[10pt]{llncs} 
\usepackage{graphicx} 
\title{Suboptimal Behavior of Bayes and MDL in Classification under Misspecification}
\usepackage{amsmath,amssymb,latexsym}

\input epsf

\newcommand{\cX}{\ensuremath{{\cal X}}}
\newcommand{\cY}{\ensuremath{{\cal Y}}}

\newcommand{\cC}{\ensuremath{{\cal C}}}
\newcommand{\cP}{\ensuremath{{\cal P}}}

\newcommand{\erateS}{\hat{e}_S}
\newcommand{\erateD}{e_D}

\newcommand{\KL}{\text{\sc KL}}

\newcommand{\ebayes}{e_D(c_{\text{\sc Bayes}(P,S)})}
\newcommand{\emap}{e_D(c_{\text{\sc map}(P,S)})}
\newcommand{\esmap}{e_D(c_{\text{\sc smap}(P,S)})}
\newcommand{\emdl}{e_D(c_{\text{\sc mdl}(P,S)})}
\newcommand{\eorb}{e_D(c_{\text{\sc orb}(P,S)})}
\newcommand{\cbayes}{c_{\text{\sc Bayes}(P,S)}}
\newcommand{\cbayesi}{c_{\text{\sc Bayes}(P,S^{i-1})}}
\newcommand{\cmap}{c_{\text{\sc map}(P,S)}}
\newcommand{\csmap}{c_{\text{\sc smap}(P,S)}}
\newcommand{\cmdl}{c_{\text{\sc mdl}(P,S)}}
\newcommand{\corb}{c_{\text{\sc orb}(P,S)}}

\newcommand{\phard}{p_{\text{hard}}}
\newcommand{\muhard}{\mu_{\text{hard}}}

\newcommand{\mhard}{m_{\text{hard}}}

\newcommand{\commentout}[1]{}
\newcommand{\journalversion}[1]{}

\author{Peter Gr\"unwald\inst{1} \and John Langford\inst{2}}
\institute{CWI Amsterdam,
\email{pdg@cwi.nl}, \texttt{www.grunwald.nl}.
\and TTI-Chicago,  \email{jcl@cs.cmu.edu}, \texttt{hunch.net/\~\/jl/}.}

\begin{document} 

\maketitle 
 
\begin{abstract} 
  We show that forms of Bayesian and MDL inference that are 
often applied to classification problems can be {\em inconsistent}.
  This means there exists a learning problem such that for all amounts
  of data the generalization errors of the MDL classifier and the
  Bayes classifier relative to the 
  Bayesian posterior both remain bounded away from the smallest
  achievable generalization error.
\end{abstract}
\section{Introduction}
Overfitting is a central concern of machine learning and
statistics. Two frequently used learning methods that in many cases
`automatically' protect against overfitting are Bayesian inference
\cite{BernardoS94} and the Minimum Description Length (MDL)
Principle \cite{Rissanen89,BarronRY98,Grunwald04}. 
We show that, when applied to
classification problems, some of the standard variations of these
two methods can be {\em inconsistent\/} in the sense that they {\em
asymptotically overfit\/}: there exist scenarios where, no matter how
much data is available, the generalization error of a classifier
based on MDL or the full Bayesian posterior does not converge to the
minimum achievable generalization error within the set of classifiers
under consideration.

\paragraph{Some Caveats and Warnings} These result must be interpreted
carefully. There exist many different versions of MDL and Bayesian
inference, only some of which are covered. For the case
of MDL, we show our result for a two-part form of MDL that has often been
used for classification, see
Section~\ref{sec:mdl}.
For the case of Bayes, our result may appear to contradict some
well-known Bayesian consistency results \cite{BlackwellD62}.
Indeed, our result only applies to a `pragmatic' use of Bayes, where
the set of hypotheses under consideration are {\em classifiers}:
functions mapping each input $X$ to a discrete class label $Y$. To
apply Bayes rule, these classifiers must be converted into conditional
probability distributions.  We do this conversion in a standard
manner, crossing a prior on classifiers with a prior on error rates
for these classifiers. This may lead to (sometimes subtly)
`misspecified' probability models not containing the `true'
distribution $D$.  Thus, our result may be restated as `Bayesian
methods for classification can be inconsistent under misspecification 
for common classification probability models'.  
The result is still interesting, since (1) even under
misspecification, Bayesian inference is known to be 
consistent under fairly broad
conditions -- we provide an explicit context in which it is not; (2) 
in practice, Bayesian inference is used frequently
for classification under misspecification -- see Section~\ref{sec:bayes}.
%
\subsection{A Preview}
\subsubsection{Classification Problems}
A classification problem is defined on an input (or feature) domain
$\cX$ and output domain (or class label) $\cY = \{0,1\}$.  The problem
is defined by a probability distribution $D$ over $\cX
\times \cY$. A classifier is a function
$c: \cX \rightarrow \cY$  The error
rate of any classifier is quantified as:
$$\erateD(c) = E_{(x,y) \sim D} I(c(x) \neq y)$$ 
where $(x,y) \sim D$ denotes
a draw from the distribution $D$ and $I(\cdot)$ is the indicator
function which is $1$ when its argument is true and $0$ otherwise.

The goal is to find a classifier which, as often as possible according to $D$,
correctly predicts the class label given the input feature.
Typically, the classification problem is solved by searching for some
classifier $c$ in a limited subset $\cC$ of all classifiers using a
sample $S = (x_1,y_1), \ldots, (x_m,y_m) \sim D^m$ 
generated by $m$ independent draws from
the distribution $D$.  Naturally, this search is guided by the
{\em empirical error rate}. This is the error rate on the subset $S$ defined by:
$$
\erateS(c) := 
E_{(x,y) \sim S} I(c(x)\neq y) = \frac{1}{m} \sum_{i=1}^m I(c(x_i)
\neq c(y_i)).
$$ where $(x,y) \sim S$ denotes a sample drawn from the uniform
distribution on $S$.  Note that $\erateS(c)$ is a random variable
dependent on a draw from $D^m$. In contrast, $\erateD(c)$ is a number (an
expectation) relative to $D$.

\subsubsection{The Basic Result}
Our basic result is that certain classifier learning algorithms may not
behave well as a function of the information they use, even when given
infinitely many samples to learn from.  The learning algorithms
we analyze are ``Bayesian classification'' (Bayes), ``Maximum a
Posteriori classification'' (MAP), and ``Minimum Description Length
classification'' (MDL).  These algorithms are
precisely defined later. Functionally they take as arguments a
training sample $S$ and a ``prior'' $P$ which is 
a
probability distribution over a set of classifiers $\cC$.
In Section~\ref{sec:result} we state our basic result,
Theorem~\ref{thm:bayesinconsistent}. The theorem has the following corollary,
indicating suboptimal behavior of Bayes and MDL:
\begin{corollary}
  \label{cor:inconsistent}
        {\bf (Classification Inconsistency) \ } There
        exists an input domain ${\cal X}$, 
a prior $P$ always nonzero on a countable set of classifiers
        $\cC$, a learning problem $D$, and a constant $K > 0$
        such that the Bayesian classifier $\cbayes$, the MAP
        classifier $\cmap$, and the MDL classifier
        $\cmdl$ are asymptotically
        $K$-suboptimal. That is, for each  ${\bf e} \in \{ \ebayes, \emap,
        \emdl \}$, we have
$$
\lim_{m \rightarrow \infty} \Pr_{S \sim D^m} \left(
 {\bf e}       
> K + \inf_{c \in \cC} \erateD(c) \right) = 1.
$$
\end{corollary}
How dramatic is this result? We may ask (1) are the priors $P$ for
which the result holds natural; (2) how large can the constant $K$
become and how small can $\inf_{c \in {\cal C}} \erateD(c)$ be? (3)
  perhaps demanding an algorithm which depends on the prior $P$ and
  the sample $S$ to be consistent (asymptotically optimal) is too strong?
The short answer to (1) and (2) is: the priors $P$ have to satisfy
  several requirements, but they correspond to priors often used in
  practice. $K$ can be quite large and $\inf_c \erateD(c)$ can be quite small
  - see Section~\ref{sec:technical} and Figure~\ref{fig:allowed}.

The answer to (3) is that there do exist simple algorithms which are
consistent. An example is the
algorithm which minimizes the Occam's Razor bound (ORB) \cite{BlumerEHW87},
Section~\ref{sec:orb}.
\begin{theorem}
  \label{thm:consistent}
        {\bf (ORB consistency) \ } For all priors $P$ nonzero on a set
        of classifiers $\cC$, for all learning problems $D$, and all constants
        $K> 0$ the ORB classifier $\corb$ 
is asymptotically $K$-optimal:
        
        $$ 
        \lim_{m \rightarrow \infty} \Pr_{S \sim D^m} \left(\eorb > K + \inf_{c \in \cC} \erateD(c) \right) = 0.$$
\end{theorem}
\commentout{
Is this relevant?
Basically any stated variant (no matter how loose) of the
Occam's Razor bound works here, which is some evidence for the
weakness of the consistency demand.  
Also, the statement
above can \emph{not} be made to work with a VC-dimension bound since
the VC-dimension of $\cC$ may be $\infty$.
}
The remainder of this paper first defines precisely what we mean by
the above classifiers. It then states the main inconsistency theorem which
implies the above corollary, as well as a theorem that provides an 
upper-bound on how badly Bayes can behave. In Section~\ref{sec:proofs}
we prove our theorems. Variations of the result are discussed in
Section~\ref{sec:technical}. A discussion of the result from a
Bayesian point of view is given in Section~\ref{sec:bayes}, and from
an MDL point of view in Section~\ref{sec:mdl}\vspace*{-0.2 cm}.
\section{Some Classification Algorithms}
\label{sec:algorithms}
The basic inconsistency result is about particular classifier learning
algorithms which we define next\vspace*{-0.2 cm}.

\subsubsection{The Bayesian Classification algorithm}
The Bayesian approach to inference starts with a prior probability
distribution $P$ over a set of distributions $\cP$ which typically
represents a measure of ``belief'' that some $p \in \cP$ is the
process generating data.  Bayes' rule states that, given sample data $S$, 
the posterior
probability $P(\cdot \mid S)$ 
that some $p$ is the process generating the data
is:
\[
P(p \mid S) =\frac{p(S)P(p)}{P(S)}.
\]
where $P(S) := E_{p \sim P} p(S)$.
In classification problems with sample size $m = |S|$, 
each $p \in \cP$ is a distribution  on $(\cX \times \cY)^n$ 
and the outcome $S=(x_1,y_1), \ldots, (x_m,y_m)$ is the sequence
of labeled examples.

If we intend to perform classification based on a set of classifiers
$\cC$ rather than distributions $\cP$, it is natural to introduce a
``prior'' $P(c)$ that a particular classifier $c: \cX \rightarrow
\{0,1\}$ is the best classifier for solving some learning problem.
This, of course, is \emph{not} a Bayesian prior in the conventional
sense because classifiers do not induce a measure over the training
data.  It is the standard method of converting a ``prior'' over
classifiers into a Bayesian prior over distributions on the
observations which our inconsistency result applies to.

One common conversion \cite{Jordan95,Tipping01,Grunwald98b}
transforms the set of classifiers $\cC$ into a simple logistic
regression model -- the precise relationship to logistic regression is
discussed in Section~\ref{sec:trafo}.  In our case $c(x) \in \{0,1\}$ is
binary valued, and then (but only then) the conversion amounts to
assuming that the error rate $\theta$ of the optimal classifier is
independent of the feature value $x$. This is known as
``homoskedasticity'' in statistics and ``label noise'' in learning
theory.  More precisely, it is assumed that, for the optimal
classifier $c \in C$, there exists some $\theta$ such that $\forall
x\,\,\, P(c(x) \neq y) = \theta.$ Given this assumption, we can
construct a conditional probability distribution $p_{c,\theta}$ over
the labels given the unlabeled data:
\begin{equation}
\label{eq:regmodel}
p_{c,\theta}(y^m  \mid x^m) =
\theta^{m \hat{e}_S(c)} (1- \theta)^{m - m \erateS(c)}.
\end{equation}
For each fixed $\theta < 0.5$, the log likelihood
$\log p_{c,\theta}(y^m \mid x^m)$ is linearly decreasing in the empirical
error that $c$ makes on $S$. 
By differentiating with respect to
$\theta$, we see that for fixed $c$, the likelihood
(\ref{eq:regmodel}) is maximized by setting $\theta :=
\erateS(c)$, giving
\begin{equation}
\label{eq:safe}
 \log \frac{1}{p_{c,\erateS(c)}(y^m \mid x^m)} = m H(\hat{e}_S(c)).
\end{equation}
where $H$ is the binary entropy $H(\mu) = - \mu \log \mu - (1- \mu)
\log (1- \mu)$, which is strictly increasing for
$\erateS(c) \in [0,0.5)$. 
We further assume that {\em some\/} distribution $p_x$ on ${\cal X}^n$
generates the $x$-values. We can apply Bayes rule to get a
posterior on $p_{c,\theta}$, denoted as $P(c,\theta \mid S)$, without knowing $p_x$, since the $p_x(x^m)$-factors cancel: 
\begin{equation}
\label{eq:bayes}
P({c,\theta} \mid S) = \frac{p_{c,\theta}(y^m | x^m)p_x(x^m)P(c,\theta)}{P(y^m \mid x^m) p_x(x^m)} =
\frac{p_{c,\theta}(y^m|x^m)P(c,\theta)}{E_{c,\theta \sim P}
  p_{c,\theta}(y^m \mid x^m) }.
\end{equation}
To make (\ref{eq:bayes})
applicable, we need to incorporate a prior measure on the joint space $\cC
\times [0,1]$ of classifiers and $\theta$-parameters. In the next
section we discuss the priors under which our theorems hold. 

Bayes rule (\ref{eq:bayes}) is formed into a classifier learning
algorithm by choosing the most likely label given the input $x$ and
the posterior $P(\cdot|S)$:
\begin{equation}
\label{eq:bayesact}
\cbayes(x) := \begin{cases}
1 & \text{if} \ E_{c,\theta \sim P(\cdot \mid S)} \/ p_{c,\theta}(Y =
1|X=x) \geq 
\frac{1}{2}, \\
0 & \text{otherwise.}
\end{cases}
\end{equation}
\subsubsection{The MAP classification algorithm}
The integrations of the full Bayesian classifier can be too
computationally intensive, so we sometimes predict using the Bayesian
Maximum A Posteriori (MAP) classifier.  This classifier is given by:
$$
\cmap = 
\arg \max_{c \in \cC} \max_{\theta \in [0,1]} P(c,\theta \mid S) = 
\arg \max_{c \in \cC} \max_{\theta \in [0,1]} 
p_{c,\theta}(y^m \mid x^m)P(c,\theta)
$$
with ties broken arbitrarily. Integration over $\theta \in [0,1]$ 
being much less problematic than summation over $c \in
\cC$, one sometimes uses a 
learning algorithm which integrates over 
$\theta$ (like full Bayes) but maximizes over $c$ (like MAP):
$$
\csmap  = 
\arg \max_{c \in \cC} P(c \mid S) =
\arg \max_{c \in \cC}
E_{\theta \sim P(\theta)} p_{c,\theta}(y^m \mid x^m) P(c \mid \theta).
$$
\subsubsection{The MDL Classification algorithm}
The MDL approach to classification is transplanted from the MDL
approach to density estimation.  There is no such thing as a
`definition' of MDL for classification because the
transplant has been performed in various ways by various authors.
Nonetheless, as we discuss in Section~\ref{sec:mdl}, most implementations are essentially equivalent to the
following algorithm \cite{QuinlanR89,Rissanen89,KearnsMNR97,Grunwald98b}:
\begin{equation}
\label{eq:cmdl}
\cmdl = 
\arg \min_{c \in \cC} \ \log \frac{1}{P(c)} + \log \binom{m}{m\hat{e}_S(c)}.
\end{equation}
The quantity minimized has a coding
interpretation: it is the number of bits required to describe the
classifier plus the number of bits required to describe the labels on
$S$ given the classifier and the unlabeled data. We call $- \log
P(c) + \log \binom{m}{m \hat{e}_S(c)}$ the {\em two-part MDL
  codelength\/} for encoding data $S$ with classifier $c$.

\section{Main Theorems} 
\label{sec:result}
\begin{figure}[tb]
\includegraphics[angle=270]{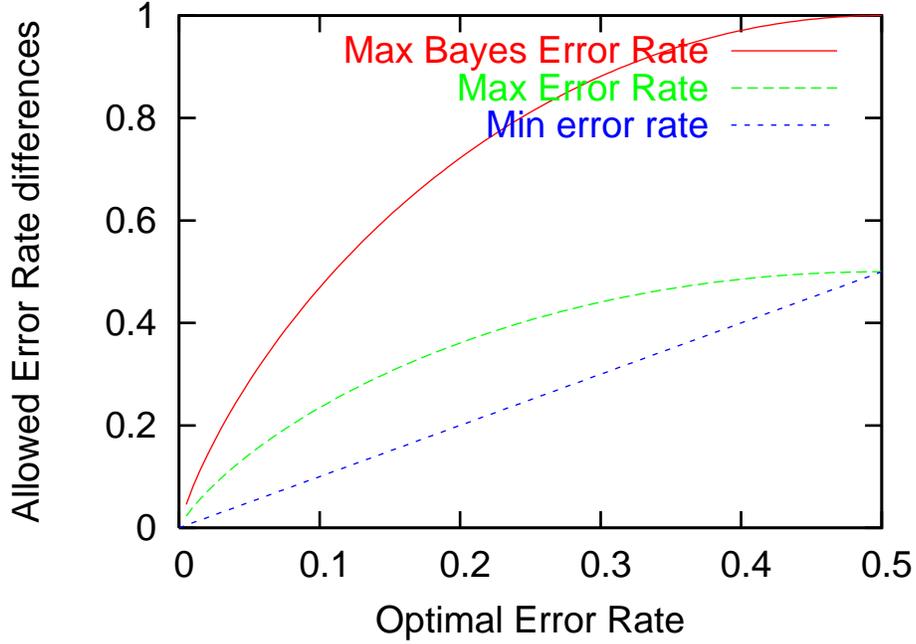}
\caption{\label{fig:allowed} A graph depicting the set of
  asymptotically allowed error rates for different classification
  algorithms. The $x$-axis depicts the optimal classifier's error rate
  $\mu$ (also shown as the straight line).  The lower curve is just
  $0.5 H(\mu)$ and the upper curve is $H(\mu)$.
Theorem~\ref{thm:bayesinconsistent} says that any $(\mu,\mu')$ 
between the straight line and
  the lower curve can be achieved for some learning
  problem $D$ and prior $P$. Theorem~\ref{thm:bayesbound} shows that
  the Bayesian learner can never have asymptotic 
error rate $\mu'$ above the upper curve.}
\end{figure}

In this section we prove the basic inconsistency theorem.  
We prove inconsistency for some countable set of classifiers
${\cal C} = \{c_0, c_1, \ldots \}$ which we define later. The
inconsistency is attained for priors with `heavy tails', satisfying
\begin{equation}
\label{eq:priora}
\log \frac{1}{P(c_k)} \leq \log k + o( \log k).
\end{equation} 
This condition is satisfied, by, for example, Rissanen's {\em
  universal prior for the integers}, \cite{Rissanen89}. 
The sensitivity of our result to the choice of prior is analyzed
further in Section~\ref{sec:technical}. 
The prior on $\theta$ can be any distribution on $[0,1]$ with a
continuously differentiable density $P$ bounded away
from $0$,
i.e. 
for some $\gamma > 0$, 
\begin{equation}
\label{eq:priorb}
\text{for all $\theta \in [0,1]$}, P(\theta) > \gamma.
\end{equation} 
For example, we may take the uniform distribution with $P(\theta)
\equiv 1$. We assume that the priors $P(\theta)$ on $[0,1]$ and the prior $P(c)$
on $\cC$ are independent, so that $P(c,\theta) = P(c) P(\theta)$.
In the theorem, $H(\mu) = - \mu \log \mu - (1- \mu) \log (1- \mu)$ 
stands for the binary entropy of a coin with bias $\mu$.
\begin{theorem}
\label{thm:bayesinconsistent}
     {\bf (Classification Inconsistency) \ } There exists an input
     space ${\cal X}$ and a countable set of classifiers ${\cal C}$
     such that the following holds: let $P$ be any
     prior satisfying (\ref{eq:priora}) and (\ref{eq:priorb}). For all $\mu
\in (0,0.5)$ and all $\mu' \in [\mu,H(\mu)/2)$,
there exists a  $D$ with
$\min_{c\in \cC} \erateD(c) = \mu$ such that, for all large $m$, all $\delta
> 0$,
\begin{eqnarray}
\Pr_{S \sim D^m} \left( \emap = \mu' \right) & \geq & 1 - a_m 
\nonumber \\
\Pr_{S \sim D^m} \left(\esmap = \mu' \right) & \geq & 1 - a_m 
\nonumber \\
\Pr_{S \sim D^m} \left(\emdl = \mu' \right) & \geq & 1 - a_m, 
\nonumber \\
\Pr_{S \sim D^m} \left(
        \ebayes \geq \mu'- \delta \right) & \geq & 1 -a_m, 
\text{\ \ \ \ where
$a_m =  3 \exp(- 2 \sqrt{m})$}.
\nonumber
\end{eqnarray}
\end{theorem}
The theorem states that Bayes is inconsistent for \emph{all large} $m$
on a fixed distribution $D$.  This is a significantly more difficult
statement than ``for all (large) $m$, there exists a learning problem
where Bayes is inconsistent''\footnote{In fact, a meta-argument can be
  made that \emph{any} nontrivial learning algorithm is `inconsistent'
  in this sense for finite $m$.}. Differentiation of $0.5 H(\mu) -\mu$
shows that the maximum discrepancy between $\emap$ and $\mu$ is
achieved for $\mu = 1/5$. With this choice of $\mu$, $0.5 H(\mu) - \mu
= 0.1609\ldots$ so that, by choosing $\mu'$ arbitrarily close to
$H(\mu)$, the discrepancy $\mu' - \mu$ comes arbitrarily close to $
0.1609\ldots$. These findings are summarized in
Figure~\ref{fig:allowed}.

How large can
the discrepancy between $\mu = \inf_c \erateD(c)$
and 
$\mu' = \erateD(\cbayes)$ be in the large $m$ limit, for general
learning problems? Our next theorem, again summarized in Figure~\ref{fig:allowed}, gives
an upperbound, namely, $\mu' < H(\mu)$:
\begin{theorem}{\bf \ (Maximal Inconsistency of Bayes)}
\label{thm:bayesbound}
Let $S^i$ be the sequence consisting of the first $i$ examples
$(x_1,y_1), \ldots, (x_i,y_i)$.  For all priors $P$ nonzero on a set
of classifiers $\cC$,
for all learning problems $D$ with $\inf_{c \in \cC} \erateD(c) = \mu$, 
for all $\delta >
0$, for all large $m$, 
with $D^m$-probability $\geq 1 - \exp(- 2 \sqrt{m})$,
$$
\frac{1}{m} \sum_{i=1}^m \bigl| y_i - \cbayesi(x_i) \bigr| \leq 
H(\mu) + \delta.
$$
\end{theorem}
The theorem says that for large $m$, the total number of mistakes when
successively classifying $y_i$ given $x_i$ made by the Bayesian
algorithm based on $S^{i-1}$, divided by $m$, is not
larger than $H(\mu)$. By the law of large numbers, it follows that
for large $m$,
$\erateD(\cbayesi(x_i))$, {\em averaged\/} over all $i$, is no larger than
$H(\mu)$. Thus, it is not ruled out that sporadically, for some
$i$, $\erateD(\cbayesi(x_i)) > H(\mu)$; but this must be `compensated'
for by most other $i$.  We did not find a proof that $\erateD(\cbayesi(x_i)) < H(\mu)$
for {\em all\/} large $i$.
\section{Proofs}
\label{sec:proofs}
In this section we present the proofs of our three theorems.
Theorem~\ref{thm:bayesinconsistent} and~\ref{thm:bayesbound} both make
use of the following lemma:
\begin{lemma}
\label{pro:basicbayes}
There exists $\gamma > 0$ such that 
for all classifiers $c$, $\alpha > 0$, $m > 0$, all $S \sim D^m$
satisfying $\alpha + 1/\sqrt{m} < \hat{e}_S(c) \leq 0.5$, all priors
satisfying (\ref{eq:priorb}): \\
\begin{multline}
\label{eq:bayesbound}
\log \frac{1}{P(y^m \mid x^m, c, \hat{e}_S(c))}
\leq  \log \frac{1}{P(y^m \mid x^m, c)} 
\leq \\
\log \frac{1}{P(y^m \mid x^m, c, \hat{e}_S(c))} + \frac{1}{2} \log {m} +
\frac{1}{2} \frac{1}{\alpha (1- \alpha)} 
- \log \gamma.
\end{multline}
\end{lemma}
\begin{proof}{\bf (sketch)\ }
For the first inequality, note
\begin{multline}
\log \frac{1}{P(y^m \mid x^m, c)} = 
\log \frac{1}{\int P(y^m \mid x^m, c,\theta) P(\theta) d\theta}
\geq 
\log \frac{1}{P(y^m \mid x^m, c, \hat{e}_S(c))}, \nonumber
\end{multline}
since the likelihood $P(y^m \mid x^m, c, \hat{e}_S(c))$ is maximized
at $\theta = \hat{e}_S(c)$. 
For the second inequality, note that
$$
\int_0^1 
P(y^m \mid x^m, c,\theta) P(\theta) d \theta 
\geq 
\int_{\hat{e}_S(c) - 1/ \sqrt{m}}^{\hat{e}_S(c) + 1 /\sqrt{m}}
\exp \bigl( \log P(y^m \mid x^m,
c,\theta) + \log P(\theta) \bigr) d \theta. 
$$
We obtain (\ref{eq:bayesbound}) by expanding $ \log P(y^m \mid
x^m, c,\theta)$ around the maximum $\theta = \hat{e}_S(c)$ using a
second-order Taylor approximation. See, \cite{BarronRY98} for further details.
\end{proof}
\subsection{Inconsistent Learning Algorithms:
Proof of Theorem~\ref{thm:bayesinconsistent}}
\label{sec:inconsistencyproof}
Below we first define the particular learning
problem that causes inconsistency. We then analyze the performance of the
algorithms on this learning problem.
\commentout{
\subsubsection{The Prior}

We use a ``Prior'' $P$ on a countable set of classifiers defined next.

Start with a measure over integers given by $Q(n) = \frac{1}{n(n+1)}$
for $n>0$.  This measure is normalized ($\sum_{n=1}^\infty
\frac{1}{n(n+1)}=1$).  Next define a measure on the classifiers given
an integer $n$ according to:
\begin{equation}
\label{eq:langfordprior}
P(c_j \mid n) = 
\begin{cases}
2^{-n+1} & \text{if} \ 
j \in \{ 2^{n-1}, 2^{n-1}+1, \ldots, 2^n - 1 \} \\
0 & \text{otherwise.}
\end{cases}
\end{equation}
This measure is also normalized ($\forall n\,\,\,\sum_{i=0}^\infty
P(c_j \mid n) = 1$), and so we can combine with $Q(n)$ to construct a
normalized prior over all classifiers.  This prior is:
\begin{equation}
\label{eq:badprior}
P(c_j) =
\begin{cases}
\frac{1}{2} & \text{if} \  j = 0\\
\frac{1}{2} \sum_{n \geq 1}^\infty P(c_j \mid n) Q(n) & \text{otherwise}
\end{cases}
\end{equation}
Note that $\sum_{j=0,1, \ldots} P(c_j) = 1$ because $Q(n)$ and
$P(c_j \mid n)$ are normalized.  The essential quality\footnote{Many
other priors also satisfy this property including Rissanen's
`universal prior for the integers' \cite{Rissanen83}, which is often
used in practical applications.} of this prior used in the proof is
that it does not decay to fast.  In particular:
\begin{equation}
\label{eq:priorreq}
\log \frac{1}{P(c_k)} \leq \log k +  o ( \log k).
\end{equation}
}
\subsubsection{The Learning Problem}
For given $\mu$ and $\mu' \geq \mu$, we 
construct a learning problem and a set of classifiers ${\cal C} =
\{ c_0, c_1, \ldots \}$ such that $c_0$ is the `good' classifier with 
$\erateD(c_0) = \mu$ and
$c_1, c_2, \ldots$ are all `bad' classifiers with $\erateD(c_j) = \mu'
\geq \mu$.
$\cX$ consists of one binary feature per classifier\footnote{ This
input space has a countably infinite size. The Bayesian posterior is
still computable for any finite \( m \) if we order the features
according to the prior of the associated classifier. We
need only consider features which have an associated prior greater
than \( \frac{1}{2^{m}} \) since the minus log-likelihood of the data
is always less than \( m \) bits. Alternatively, we could use
stochastic classifiers and a very small input space.}, and
the classifiers simply output the value of their special feature.
The underlying distribution $D$ is constructed in terms of 
$\mu$ 
and $\mu'$
and a proof parameter
$\muhard \geq \frac{1}{2}$ (the error rate for ``hard'' examples).
To construct an example $(x,y)$, we first flip a fair coin to
determine $y$, so $y =1$ with probability $1/2$. We then flip a coin with bias
$\phard :=\frac{\mu'}{\muhard}$ which determines if this is a ``hard''
example or an ``easy'' example.  Based upon these two coin flips, each $x_j$ is independently generated based on the following $3$ cases.
\begin{enumerate}
\item For a ``hard'' example, and for each classifier $c_j$ 
with $j \geq 1$, set $x_j = |1 - y|$ with probability $\muhard$ and
$x_j = y$  otherwise.
\item For an ``easy'' example, and every $j \geq 1$ set $x_j = y$.
\item For the ``good'' classifier $c_0$ (with true error rate $\mu$),
  set $x_0 = |1- y|$ with probability $\mu$ and $x_0 = y$ otherwise.
\end{enumerate}
The error rates of each classifier are $\erateD(c_0) = \mu$ and
$\erateD(c_j) = \mu'$ for all $j \geq 1$.
\subsubsection{Bayes and MDL are inconsistent}
We now prove Theorem~\ref{thm:bayesinconsistent}.
In Stage 1 we show that there exists a $k_m$ such that for every
value of $m$, with probability converging to 1, there exists
some `bad' classifier $c_j$ with $0 < j \leq k_m$ that has $0$
empirical error. In Stage 2 we show that the prior of this classifier
is large enough so that its posterior is exponentially larger than
that of the good classifier $c_0$, showing the convergence $\emap
\rightarrow \mu'$. In Stage 3 we sketch the convergences $\esmap
\rightarrow \mu', \emdl \rightarrow \mu', \ebayes \rightarrow \mu'$.
\paragraph{Stage 1}
Let $\mhard$ denote the number of hard examples generated within a
sample $S$ of size $m$. Let $k$ be a positive integer and $\cC_k =
\{c_j \in \cC: 1 \leq j \leq k  \}$.  For all $\epsilon > 0$
and $m \geq 0$, we have:
\begin{eqnarray}
\lefteqn{\Pr_{S \sim D^m}(\forall c \in \cC_k: \ \erateS(c) > 0)} & &
\nonumber  \\ 
& \overset{(a)}{=} & 
\Pr_{S \sim D^m} \left( \forall c \in \cC_k: \ \erateS(c) > 0 
\mid \frac{\mhard}{m} > \phard + \epsilon \right) 
\Pr_{S\sim D^m} \left( \frac{\mhard}{m} > \phard + \epsilon \right) 
\nonumber \\
& + & \Pr_{S \sim D^m} \left( \forall c \in \cC_k: \ \hat{e}_S(c) > 0 
\mid \frac{\mhard}{m} \leq \phard + \epsilon \right) 
\Pr_{S \sim D^m} \left( \frac{\mhard}{m} \leq  \phard + \epsilon
\right) 
\nonumber  \\
& \overset{(b)}{\leq} & e^{-2m \epsilon^2} + 
\Pr_{S \sim D^m}\left( \forall c \in \cC_k: \ \erateS(c) > 0 \mid
  \frac{\mhard}{m} \leq \phard + \epsilon \right) 
\nonumber \\
\label{eq:monday}
& \overset{(c)}{\leq} &  e^{-2m \epsilon^2} + 
(1 - (1 - \muhard)^{m \left( \phard + \epsilon \right)})^k 
\overset{(d)}{\leq} 
e^{-2m \epsilon^2} + e^{-k(1 - \muhard)^{m\left( \phard +
      \epsilon \right)}}.
\end{eqnarray}
Here (a) follows because $P(a)= \sum_b P(a|b)P(b)$.
(b) follows by $\forall a,P:\,\,P(a) \leq 1$ and the Chernoff bound.
(c) holds
since $(1 - (1 - \muhard)^{m(\phard + \epsilon)})^k$ is monotonic in
$\epsilon$, and (d) by 
$\forall x \in [0,1], k>0:\,\,\, (1-x)^k \leq e^{-kx}$.
We now set $\epsilon_m := m^{-0.25}$ and 
$k(m) = \frac{2m \epsilon_m^2}{(1 -
\muhard)^{m\left( \phard + \epsilon_m \right)}}$. Then
(\ref{eq:monday}) becomes
\begin{equation}
\label{eq:friday}
\Pr_{S \sim D^m}(\forall c \in \cC_{k(m)}: 
\ \erateS(c) > 0) \leq  2 e^{-2\sqrt{m}}
\end{equation}
On the other hand, by the Chernoff bound we have
$
\Pr_{S \sim D^m}(\erateS(c_0) < \erateD(c_0) - \epsilon_m) \leq 
e^{-2 \sqrt{m}} 
$
for the optimal classifier $c_0$. Combining this
with (\ref{eq:friday}) using the union bound, we get that, with
$D^m$-probability larger than $1- 3 e^{-2\sqrt{m}}$, the following
event holds:
\begin{equation}
\label{eq:saturday}
\exists c \in \cC_{k(m)}: 
\ \erateS(c) = 0 \text{\ and \ } \erateS(c_0) \geq \erateD(c_0) -
\epsilon_m.
\end{equation}
\paragraph{Stage 2}
In the following derivation, we assume  that the large probability event
(\ref{eq:saturday}) holds. We show that this implies that  for large $m$, the posterior on
some $c^* \in \cC_{k(m)}$ with $\erateS(c^*) = 0$ is greater than the
posterior on $c_0$, which implies that the MAP algorithm is inconsistent.
Taking the log of the posterior ratios, we get:
\begin{multline}
\label{eq:odds}
\ \hspace*{-0.6 cm}
\log \frac{\max_{\theta}
P(c_0, \theta \mid x^m,y^m)}{\max_{\theta} P(c^*, \theta \mid x^m,y^m)} = 
\log  
\frac{\max_{\theta} P(c_0) P(\theta) P(y^m \mid x^m, c_0,  \theta)}
{\max_{\theta} P(c^*) P(\theta) P(y^m \mid x^m, c^*, \theta)}
= \\
\ \hspace*{-0.3 cm} \log \max_{\theta} P(c_0) P(\theta) P(y^m \mid x^m, c_0,  \theta)
-
\log \max_{\theta} P(c^*) P(\theta) P(y^m \mid x^m, c^*, \theta).
\end{multline}
Using (\ref{eq:safe}) we see that
the leftmost term is no larger than
\begin{multline}
\label{eq:leftie}
\log \ \bigl( \max_{\theta} P(c_0) P(\theta) \bigr) \cdot 
\bigl(\max_{\theta'}  P(y^m \mid x^m, c_0, {\theta'}) \bigr)
= - m H(\erateS(c_0)) + O(1) \leq \\
-m H(\erateD(c_0)) - K m \epsilon_m + O(1)
= - m H(\mu) - {m}^{0.75} K  + O(1) 
\end{multline}
where $K$ is some constant. 
The last line follows because $H(\mu)$ is continuously differentiable 
in a small enough neighborhood around $\mu$.

For the rightmost term in (\ref{eq:odds}), by the condition on prior $p(\theta)$, (\ref{eq:priorb}), 
\begin{equation}
\label{eq:rightie}
-
\log \max_{\theta} P(c^*) P(\theta) P(y^m \mid x^m, c^*, \theta) \leq
- \log P(c^*) + \log \gamma.
\end{equation}
Using
condition (\ref{eq:priora}) on  prior $P(c^*)$ and using
$c^* \in \cC_{k(m)}$, we find:
\begin{equation}
\label{eq:priorc}
\log \frac{1}{P(c^*)} \leq \log k(m) + o ( \log k(m)),
\end{equation}
where $\log k(m) = \log {2}{\sqrt{m}} - 
(m \phard + m^{0.25}) \log (1 - \muhard)$. 
Choosing $\muhard = 1/2$, this becomes
$\log k(m) =  \frac{1}{2} \log m  + 2 m \mu'  + m^{0.25} + O(1)$.
Combining this with (\ref{eq:priorc}), we find that 
\begin{equation}
\label{eq:key}
 \log
\frac{1}{P(c^*)} \leq 2 m \mu' + o(m) 
\end{equation}
which implies that 
(\ref{eq:rightie}),
is no larger than $ 2 m \mu' + o(m)$.  Since $\mu' < H(\mu)/2$, the
difference between the leftmost term (\ref{eq:leftie})  and the 
rightmost term (\ref{eq:rightie}) in (\ref{eq:odds}) is
less than $0$ for large $m$, implying that then $\emap = \mu'$. We
derived all this from (\ref{eq:saturday}) which holds with
probability $\geq 1 - 3 \exp(- 2 \sqrt{m})$. Thus, 
for all large $m$,
$ \Pr_{S \sim D^m} \left( \cmap = \mu' \right) \geq 1-
3 \exp(- 2 \sqrt{m}),$
and the result follows.
\paragraph{Stage 3} {\bf (sketch)}
The proof that the integrated MAP classifier $\csmap$ is inconsistent
is similar to the proof for $\cmap$ that we just gave, except that
(\ref{eq:odds}) now becomes
\begin{equation}
\label{eq:smapodds}
\log P(c_0) P(y^m \mid x^m, c_0)
-
\log P(c^*) P(y^m \mid x^m, c^*).
\end{equation}
By Lemma~\ref{pro:basicbayes} we see that, if (\ref{eq:saturday})
holds, the difference between (\ref{eq:odds}) and (\ref{eq:smapodds})
is of order
$O(\log m)$. The proof then proceeds exactly as for the MAP case.

To prove inconsistency of $\cmdl$, note that the MDL code length of $y^m$ given $x^m$ according to $c_0$ is given by
$\log \binom{m}{m \hat{e}_S(c_0)}$. If (\ref{eq:saturday}) holds, then a simple Stirling's approximation as in \cite{Grunwald98b} or \cite{KearnsMNR97} shows that $\log \binom{m}{m \hat{e}_S(c_0)} = m H(\hat{e}_S(c_0)) - O(\log m)$. Thus, the difference between two-part codelengths achieved by $c_0$ and $c^*$ is given by
\begin{equation}
\label{eq:mdlodds}
- m H(\hat{e}_{S}(c_0)) + O(\log m) - \log P(c^*).
\end{equation} 
The proof then proceeds as for the MAP case, with (\ref{eq:odds})
replaced by (\ref{eq:mdlodds}) and a few immediate adjustments.

To prove inconsistency of $\cbayes$, we take $\muhard$ not equal to
$1/2$ but to $1/2 + \delta$ for some small $\delta > 0$. By taking
$\delta$ small enough, the proof for $\cmap$ above goes through
unchanged so that, with probability  $\geq 1 -
3 \exp(-2 \sqrt{m})$, the Bayesian
posterior puts all its weight, except for an exponentially small part,
on a mixture of distributions $p_{c_j,\theta}$ 
whose Bayes classifier has 
error rate $\mu'$ and error rate on hard examples $> 1/2$. It can
be shown that this implies that
for large $m$, the classification error $\cbayes$ converges to
$\mu'$; we omit details.   
\commentout{
What we have left to prove is the inconsistency of the full Bayes
classifier and in particular the surprising fact that the error rate
of Bayes can double.

We need a different prior in order to prove this statement.  If the
prior has the property that the ratio between the size of the set of
classifiers at one $n+1$ and at size $n$ grows by a factor of $2$ with
each increment, then we should be able to prove something like this.
In particular we should be able to show that it is infinitely often
the case that for some critical $n$ there are:

1) No zero error classifiers for smaller $n$.

2) A growing-with-$n$ number of zero error classifiers in the set
   covered by $n$.

Then, we should be able to apply concentration inequalities to the set
of zero error classifiers at the fixed $n$ implying that the full
bayes classifier errs on all hard examples.  Getting this right
seems mathematically tricky.  
}
\subsection{A Consistent Algorithm: Proof of Theorem~\ref{thm:consistent}}
\label{sec:orb}
In order to prove the theorem, we first state
the Occam's Razor Bound classification algorithm, based on minimizing the bound given by the following theorem.
\begin{theorem} (Occam's Razor Bound) {\rm  \cite{BlumerEHW87}} 
\label{thm:orb}
For all priors $P$ on a countable set of classifiers $\cC$, 
for all distributions $D$, with
probability $1-\delta$:
$$ \forall c:\,\,\,\,\erateD(c) \leq \erateS(c) + \sqrt{\frac{\ln \frac{1}{P(c)} + \ln \frac{1}{\delta}}{2m}}.$$
\end{theorem}
We  state the algorithm here in a suboptimal form, which good enough for our purposes (see \cite{McAllester99} for more sophisticated versions): 
$$ \corb := \arg \min_{c \in \cC} \erateS(c) + \sqrt{\frac{\ln
    \frac{1}{P(c)} + 
\ln {m}}{2m}}.$$
\paragraph{Proof of Theorem~\ref{thm:consistent}}
Set $\delta_m := 1/m$. It is easy to see that 
$$
\min_{c \in \cC} \erateD(c) + \sqrt{\frac{\ln
    \frac{1}{P(c)} + 
\ln {m}}{2m}}
$$
is achieved for at least one $c \in \cC = \{ c_0, c_1, \ldots
\}$. Among all $c_j \in \cC$ achieving the minimum, let 
$\tilde{c}_{m}$ be the one with smallest index $j$. By the Chernoff
bound, we have with probability at least $1-
\delta_m = 1- 1/m$,
\begin{equation}
\label{eq:chernoccam}
\erateD(\tilde{c}_{m}) \geq \erateS(\tilde{c}_{m}) -
\sqrt{\frac{\ln (1/\delta_m)}{2m}} = 
\erateS(\tilde{c}_{m}) -
\sqrt{\frac{\ln m}{2m}},
\end{equation}
whereas by Theorem~\ref{thm:orb}, with probability at least
$1 - \delta_m = 1- 1/m$,
$$
\erateD(\corb) \leq \min_{c \in \cC} 
\erateS(c) + \sqrt{\frac{- \ln
    {P(c)} + 
\ln {m}}{2m}} \leq \erateS(\tilde{c}_{m}) +
\sqrt{\frac{- \ln {P(\tilde{c}_{m})} + 
\ln {m}}{2m}}.
$$
Combining this with (\ref{eq:chernoccam}) using the union bound, we find that 
$$\erateD(\corb) \leq \erateD(\tilde{c}_{m}) +  \sqrt{\frac{- \ln
    {P(\tilde{c}_{m})} + 
\ln {m}}{2m}} + \sqrt{\frac{\ln m}{2m}},$$
with probability at least $1 - 2/m$.
\enlargethispage{\baselineskip}
The theorem follows upon noting that the right-hand side of this
expression converges to $\inf_{c \in \cC} \erateD(c)$ with increasing $m$\vspace*{-0.3 cm}. 
\commentout{Pick any $\epsilon$.  Pick any classifier $c'$ with $P(c') > 0$ satisfying
$\erateD(c) \leq \epsilon/2 + \inf_{c \in \cC} \erateD(c)$.  
For any fixed
$m$ we have:
\begin{multline}
\Pr_{S \sim D^m} (\eorb > \inf_{c \in \cC} \erateD(c))  \leq \\ 
\Pr_{S \sim D^m} \left( \exists c \in \cC:\,\,\erateD(c) 
> \erateS(c) + \sqrt{\frac{\ln \frac{1}{P(c)} 
+ \ln \frac{m}{\delta}}{2m}}\right) \\
+ \Pr_{S \sim D^m} \left( \forall c:\erateD(c) \leq \epsilon/2 + 
\inf_{c \in \cC} \erateD(c):\,\,\epsilon + \inf_{c \in \cC} \erateD(c)
< \erateS(c) + \sqrt{\frac{\ln \frac{1}{P(c)} 
+ \ln \frac{m}{\delta}}{2m}}\right) \nonumber \end{multline}
Since $\eorb$ is works by definition when the bound holds (negation of first term) and quantity minimized is less than 
$\epsilon + \inf_{c\in \cC}$ (second term).  

In the second term, we can loosen the inequality by choosing some
particular $c$ with $ \erateD(c) \leq \epsilon/2 + \inf_{c \in \cC}$.
Doing this, we get:
\begin{multline}
\leq \frac{\delta}{m} + \Pr_{S \sim D^m}\left( \erateS(c) -\erateD(c)
  > \epsilon/2 - \sqrt{\frac{\ln \frac{1}{P(c)} + \ln
      \frac{m}{\delta}}{2m}} \right) \leq \\
\frac{\delta}{m} + \exp \left(-2m \left (\epsilon/2 -
    \sqrt{\frac{\ln \frac{1}{P(c)} + \ln
        \frac{m}{\delta}}{2m}}\right)^2 \right) \nonumber
\end{multline}
where the first inequality follows from the Occam's Razor bound.  
and the second term from the Chernoff bound.  The limit as
$m\rightarrow \infty$ of the above is $0$.
}
\subsection{Proof of Theorem~\ref{thm:bayesbound}}
\enlargethispage{\baselineskip}
Without loss of generality assume that $c_0$ achieves $\min_{c \in \cC} \erateD(c)$. Consider both the $0/1$-loss and the log loss of sequentially
  predicting with the Bayes predictive distribution $P(Y_i =  \cdot
  \mid X_i = \cdot , S^{i-1})$ given by
$P(y_i  \mid  x_i, S^{i-1}) =  
E_{c,\theta \sim P(\cdot \mid S^{i-1})} \/ p_{c,\theta}(y_i|x_i).$
Every time $i \in
  \{1, \ldots, m \}$ that the Bayes classifier based on $S^{i-1}$
  classifies $y_i$ incorrectly,
  $P(y_i \mid x_i, S^{i-1})$ must be $\leq 1/2$ so that 
$ - \log P(y_i \mid x_i, S^{i-1}) \geq 1$.
Therefore\vspace*{-0 mm}, 
\begin{equation}
\label{eq:champagne}
\sum_{i=1}^m  - \log P(y_i \mid x_i,  S^{i-1}) \geq 
\sum_{i=1}^m | y_i - \cbayesi(x_i)|.
\end{equation}
On the other hand we have\vspace*{-0.2 mm}\ 
\begin{multline}
\label{eq:coffee}
\sum_{i=1}^m  -\log P(y_i \mid x_i, S^{i-1}) = 
 -\log \prod_{i=1}^m P(y_i \mid x_i, x^{i-1},y^{i-1}) =  \\
 -\log \prod_{i=1}^m P(y_i \mid x^m, y^{i-1}) = 
- \log \prod_{i=1}^m \frac{P(y^i | x^m)}{P(y^{i-1}| x^{m})} {=}
  - \log P(y^m \mid x^m) = \\
 - \log \sum_{j = 0,1,2 \ldots} P(y^m \mid x^m, c_j)
P(c_j)   {\leq} 
 - \log P(y^m \mid x^m, c_0) - \log P(c_0),
\end{multline}
where the inequality follows because a sum is larger than each of its
terms. By the Chernoff bound, for all small enough $\epsilon > 0$,
with probability larger than $1 - 2 \exp(-2 m \epsilon^2)$, we have $|
\erateS(c_0) - \erateD(c_0) | < \epsilon$. We now set $\epsilon_m =
m^{-0.25}$. Then, using Lemma~\ref{pro:basicbayes}, with probability
larger than $1- 2 \exp(- 2 \sqrt{m})$, for all large $m$
(\ref{eq:coffee}) is less than or equal to\vspace*{-0.3 cm}\ 
\begin{multline}
\label{eq:tea}
 - \log P(y^m \mid x^m, c_0, \hat{e}(c_0)) +
\frac{1}{2} \log m + C_{m} \overset{(a)}{=} 
m H(\erateS(c_0)) + \frac{1}{2} \log m + C_{m} \leq \\
m H(\erateD(c_0)) + Km^{0.75}  + \frac{1}{2} \log m + C_{m},
\end{multline}
where $C_{m} = (\erateD(c_0) - \epsilon_m - m^{- 0.5})^{-1}(1-
\erateD(c_0) + \epsilon_m + m^{- 0.5})^{-1}$ and $K$ is a constant not
depending on $S= S^m$.  Here $(a)$ follows from Equation~\ref{eq:safe}
and (b) follows because $H(\mu)$ is continuously differentiable in a
neighborhood of $\mu$.

\enlargethispage{\baselineskip}
Combining (\ref{eq:tea}) with (\ref{eq:champagne}) and using $C_m = O(1)$
we find  that with probability $\geq 1 - \exp(- 2 \sqrt{m})$,
$
\sum_{i=1}^m | y_i - \cbayesi(x_i)| \leq 
m H(\erateD(c_0)) + o(m),
$ QED\vspace*{-0.3 cm}.
\section{Technical Discussion} 
\subsection{Variations of Theorem~\ref{thm:bayesinconsistent} and 
dependency on the prior}
\label{sec:technical}
\enlargethispage{\baselineskip}
\paragraph{Prior on classifiers}
The requirement (\ref{eq:priora}) that
$- \log P(c_k) \geq \log k +  o ( \log k)$
is needed to obtain (\ref{eq:key}), which is the key inequality in the proof of Theorem~\ref{thm:bayesinconsistent}.  If $P(c_k)$ decreases at polynomial rate, 
but at a degree $d$ larger than one, i.e. if 
\begin{equation}
\label{eq:dtail}
-\log P(c_k) = d \log k +  o ( \log k),
\end{equation} 
then a variation of Theorem~\ref{thm:bayesinconsistent} still applies
but the maximum possible discrepancies between $\mu$ and $\mu'$ become
much smaller: essentially, if we require $\mu \leq \mu' < \frac{1}{2d}
H(\mu)$ rather than $\mu \leq \mu' < \frac{1}{2} H(\mu)$ as in
Theorem~\ref{thm:bayesinconsistent}, then the argument works for all
priors satisfying (\ref{eq:dtail}). Since the derivative $d H(\mu) /
d\mu \rightarrow \infty$ as $\mu \downarrow 0$, by setting $\mu$ close
enough to $0$ it is possible to obtain inconsistency for any fixed
polynomial degree of decrease $d$. However, the higher $d$, the smaller $\mu =
\inf_{c \in \cC} \erateD(c)$ must be to get any inconsistency with our
argument\vspace*{-1 mm}. 
%
\paragraph{Prior on error rates}
Condition (\ref{eq:priorb}) on the prior on the error rates is
satisfied for most reasonable
priors. 
Some approaches to applying MDL to classification problems amount to
assuming priors of the form $p(\theta^*) = 1$ for a single $\theta^*
\in [0,1]$ (Section~\ref{sec:mdl}). In that case, we can still
prove a version of Theorem~\ref{thm:bayesinconsistent}, but the
maximum discrepancy between $\mu$ and $\mu'$ may now be either larger
or smaller than $H(\mu)/2 - \mu$, depending on the choice of
$\theta^*$.
\subsection{Properties of the transformation 
  from classifiers to distributions} 
\label{sec:trafo}
\paragraph{Optimality and Reliability}
Assume that the conditional distribution of $y$ given $x$ according to the 
`true' underlying distribution $D$ is defined for all $x \in \cX$, and let $p_D(y|x)$ denote its mass function. 
Define $\Delta(p_{c,\theta})$ as the Kullback-Leibler (KL) divergence \cite{CoverT91} between 
$p_{c,\theta}$ and the `true' conditional distribution $p_D$:
$$
\Delta(p_{c,\theta}) := 
 \KL(p_D \| p_{c,\theta}) =
E_{(x,y) \sim D} [ - \log p_{c,\theta}(y|x) + \log p_{D}(y|x)]. 
$$
\begin{proposition}
\label{pro:trafo}
Let $\cC$ be any set of classifiers, and let $c^*
\in \cC$ achieve \\ $\min_{c \in \cC} \erateD(c)$ $= \erateD(c^*)$.
\begin{enumerate}
\item If $\erateD(c^*) < 1/2$, then 
$$\min_{c,\theta} \Delta(p_{c,\theta}) \text{\ is uniquely achieved for 
\ } (c,\theta) = (c^*,\erateD(c^*)).$$
\item $\min_{c,\theta} \Delta(p_{c,\theta}) = 0$ iff $p_{c^*,\erateD(c^*)}$ is `true', i.e. if
$\forall x, y: \/ p_{c^*,\erateD(c^*)}(y|x) = p_D(y|x)$.
\end{enumerate}
\end{proposition}
Property 1 follows since for each fixed $c$, $\min_{\theta \in
  [0,1]} \Delta(p_{c,\theta})$ is uniquely achieved for $\theta =
\erateD(c)$ (this follows by differentiation) and satisfies $
\min_{\theta} \Delta(p_{c,\theta}) = \Delta(p_{c,\erateD(c)}) =
H(\erateD(c)) - K_D, $ where $K_D = E [\log p_{D}(y|x)]$ does not
depend on $c$ or $\theta$, and $H(\mu)$ is monotonically increasing for $\mu < 1/2$. Property 2 follows from the information
inequality \cite{CoverT91}.

Proposition~\ref{pro:trafo} implies that our transformation is a good
candidate for turning classifiers into probability distributions.

Namely, let $\cP = \{ p_\alpha : \alpha \in A \}$ be a set of i.i.d.
distributions indexed by parameter set $A$ and let $P(\alpha)$ be a
prior on $A$.  By the law of large numbers, for each $\alpha \in A$,
$m^{-1} \log p_{\alpha}(y^m \mid x^m)P(\alpha) \rightarrow \KL(p_D \|
p_{\alpha})$. By Bayes rule, this implies that if the class $\cP$ is
`small' enough so that the law of large numbers holds {\em uniformly\/} 
for all $p_{\alpha} \in
\cP$, then for all $\epsilon > 0$, the Bayesian posterior will
concentrate, with probability 1, on the set of distributions in $\cP$ within
$\epsilon$ of the $p^* \in \cP$ minimizing KL-divergence to $D$. In
our case, if $\cC$ is `simple' enough so that the corresponding $\cP =
\{ p_{c,\theta} : c \in \cC, \theta \in [0,1] \}$ admits uniform
convergence \cite{Grunwald98b}, then 
%
the Bayesian posterior asymptotically concentrates  on 
the $p_{c^*,\theta^*} \in \cP = \{ p_{c,\theta} \}$ closest to $D$ in
KL-divergence. By Proposition~\ref{pro:trafo}, this $p_{c^*,\theta^*}$
corresponds to the $c^* \in \cC$ with smallest generalization error
rate $\erateD(c^*)$ ($p_{c^*,\theta^*}$ is {\em optimal\/} for
$0/1$-loss), and for the $\theta^* \in [0,1]$ with $\theta^* =
\erateD(c^*)$ ($p_{c^*,\theta^*}$ gives a {\em reliable\/} impression
of its prediction quality). This convergence to an optimal and
reliable $p_{c^*,\theta^*}$ will happen if, for example, $\cC$ has
finite VC-dimension \cite{Grunwald98b}. We can only get trouble as in
Theorem~\ref{thm:bayesinconsistent} if we allow $\cC$ to be of
infinite VC-dimension.
  
\paragraph{Analogy to Regression} In ordinary (real-valued) 
regression, ${\cal Y } = {\mathbb R}$,
and one tries to learn a function $f \in {\cal F}$ from the data. Here ${\cal F}$ is 
a set of candidate functions
${\cal X} \rightarrow {\cal Y}$. In order to apply Bayesian inference
to this problem, one assumes a probability model $\cP$ expressing 
$Y = f(X) + Z$, where
$Z$ is independent noise with mean $0$ and variance $\sigma^2$. $\cP$ then 
consists of conditional density functions $p_{f,\sigma^2}$, one for each $f \in {\cal F}$ and $\sigma^2 > 0$. It is
well known that if one assumes $Z$ to be normally distributed
independently of $X$, then the $p_{f,\sigma^2}$ become Gaussian densities and 
the likelihood becomes {\em a linear function of
  the mean squared error\/} \cite{Rissanen89}:
\begin{equation}
\label{eq:gauss}
- \ln p_{f,\sigma^2}(y^n \mid x^n) = \beta_{\sigma} \sum_{i=1}^n (y_i - f(x_i))^2 + 
n \ln Z(\beta_{\sigma}).
\end{equation}
where we wrote $\beta_{\sigma} = 1/2 \sigma^2$ and 
$Z(\beta) = \int_{y \in \cY}
\exp(- \beta y^2)d y$. Because least squares is an intuitive,
mathematically well-behaved and easy to perform procedure, it is often
assumed in Bayesian regression that the noise is normally distributed
-- even in cases where in reality, it is not
\cite{Grunwald98b,KleijnvdV04}.

Completely analogously to the Gaussian case, our transformation maps classifiers $c$ and noise rates $\theta$ to distributions $p_{c,\theta}$ so that the likelihood becomes a {\em linear function of the
  $0/1$-error}, since it can be written as:
\begin{equation}
\label{eq:entrop}
- \ln p_{c,\theta}(y^n \mid x^n) = \beta_{\theta} 
\sum_{i=1}^n |y_i - c(x_i)| + 
n \ln Z(\beta_{\theta}).
\end{equation}
where we wrote $\beta_\theta = \ln (1- \theta) - \ln \theta$ and $Z(\beta) =
\sum_{y \in \cY} \exp(- \beta y)$ \cite{Grunwald98b,MeirM95}.  Indeed, the
models $\{ p_{c,\theta} \}$ are a special case of {\em
  logistic regression models}, which we now define:
\paragraph{Logistic regression interpretation}
let $\cC$ be a
set of functions $\cX \rightarrow \cY$, where $\cY \subseteq {\mathbb
  R}$ ($\cY$ does not need to be binary-valued). The corresponding
logistic regression model is the set of conditional distributions $\{
p_{c,\beta} : c \in {\cal C} ; \beta \in {\mathbb R} \}$ of the
form
\begin{equation}
\label{eq:logit}
p_{c,\beta}(1 \mid x) := \frac{e^{-\beta c(x)}}{1 + e^{-\beta c(x)}} \ \ ; \ \ 
p_{c,\beta}(0 \mid x) := \frac{1}{1 + e^{-\beta c(x)}}.
\end{equation}
This is the standard construction used to convert classifiers with
real-valued output such as support vector machines and neural networks
into conditional distributions \cite{Jordan95,Tipping01}, 
so that Bayesian
inference can be applied.  By setting $\cC$ to be a set of
$\{0,1\}$-valued classifiers, and substituting $\beta = \ln (1-
\theta) - \ln \theta$ as in (\ref{eq:entrop}), we see that our
construction is a special case of the logistic regression
transformation (\ref{eq:logit}). It may seem that (\ref{eq:logit})
does not treat $y = 1$ and $y= 0$ on equal footing, but this is not
so: we can alternatively define a symmetric version of
(\ref{eq:logit}) by defining, for each $c \in \cC$, a corresponding $c': \cX
\rightarrow \{-1,1\}$, $c'(x) := 2c(x) -1$. Then we can set
\begin{equation}
\label{eq:logitb}
p_{c,\beta}(1 \mid x) := \frac{e^{-\beta c(x)}}{e^{\beta c(x)} + e^{-\beta c(x)}} \ \ ; \ \ 
p_{c,\beta}(-1 \mid x) := \frac{e^{\beta c(x)}}{e^{\beta c(x)} + e^{-\beta c(x)}}.
\end{equation}
By setting $\beta' = 2 \beta$ we see that $p_{c,\beta}$ as in
(\ref{eq:logit}) is identical to $p_{c,\beta'}$ as in
(\ref{eq:logitb}), so that the two models really coincide.  
%
\section{Interpretation from a Bayesian perspective}
\label{sec:bayes}
\subsubsection{Bayesian Consistency}
It is well-known that Bayesian inference is strongly consistent under very
broad conditions. For example, when applied to our setting, 
the celebrated Blackwell-Dubins
consistency theorem \cite{BlackwellD62} says the following. Let $\cC$
be countable and suppose
$D$ is such that, for some $c^* \in \cC$ and $\theta^* \in [0,1]$,
$p_{c^*,\theta^*}$ is equal to $p_D$, the true distribution/ mass function of $y$ given
$x$. 
Then with $D$-probability
1,  the Bayesian posterior concentrates on $c^*$:
$ \lim_{m \rightarrow \infty} 
P( c^* \mid S^m) = 1$.

Consider now the learning problem underlying
Theorem~\ref{thm:bayesinconsistent} as described in
Section~\ref{sec:inconsistencyproof}. Since $c_0$ achieves $\min_{c
  \in \cC} \erateD(c)$, it follows by part 1 of
Proposition~\ref{pro:trafo} that $\min_{c,\theta}
\Delta(p_{c,\theta}) = \Delta(p_{c_0,\erateD(c_0)})$. {\em If\/} $
\Delta(p_{c_0,\erateD(c_0)})$ were $0$, then by part 2 
of Proposition~\ref{pro:trafo}, 
Blackwell-Dubins would apply, and we would have $P(c_0
\mid S^m) \rightarrow 1$. Theorem~\ref{thm:bayesinconsistent} states
that this does {\em not\/} happen.  It follows that the premisse
$\Delta(p_{c_0,\erateD(c_0)}) = 0$ must be false. But since 
$\Delta(p_{c,\theta})$ is minimized for $(c_0,\erateD(c_0))$, the 
Proposition implies that for {\em no\/} $c \in \cC$ and {\em no\/}
$\theta \in [0,1]$, $p_{c,\theta}$ is equal to $p_D(\cdot|\cdot)$ - in statistical terms, the model $\cP = \{ p_{c,\theta} : c \in \cC, \theta \in [0,1] \}$ is {\em misspecified}. Thus, our result can be interpreted in two
ways:
\begin{enumerate}
\item {\em `ordinary' Bayesian inference can be inconsistent under
    misspecification\/}: We exhibit a simple logistic regression model
  $\cP$ and a true
  distribution $D$ such that, with probability 1, the Bayesian
  posterior does not converge to the distribution
  $p_{c_0,\erateD(c_0)} \in \cP$ that minimizes, among all $p \in
  \cP$, the KL-divergence to $D$, even though $p_{{c}_0,\erateD(c_0)}$
  has substantial prior mass and is {\em partially\/} correct in the
  sense that $c_0$, the Bayes optimal classifier relative to
  $p_{{c}_0,\erateD(c_0)}$, has true error rate $\erateD(c_0)$, which
  is the {\em same\/} true error rate that it would have if
  $p_{{c}_0,\erateD(c_0)}$ were `true'.
\item {\em `pragmatic' Bayesian inference for classification can be
    suboptimal\/}: a standard way to turn classifiers into
  distributions so as to make application of Bayesian
  inference possible may give rise to suboptimal performance\vspace*{-0.0 cm}. 
\end{enumerate}
\subsubsection{Two types of misspecification} $p_{c_0,\erateD(c_0)}$ can be misspecified in
  {\em two different ways}. $p_{c_0,\erateD(c_0)}$ expresses that $y =
  c_0(x) \ {\tt xor} \  z$ where $z$ is a noise bit generated 
independently of $x$. 
This statement may be wrong either because (a)  $c_0$ is not the Bayes optimal
classifier according to $D$; or (b) $c_0$ is Bayes optimal, but $z$ is
  dependent on $x$ under $D$.  The way we defined our learning problem $D$
(Section~\ref{sec:inconsistencyproof}) is an example of case (a). But
we could have equally defined $c_0$ as follows: we replace step 3 of
the generation of input values $x_j$ by the following procedure: for
an easy example, we set $x_0 = y_0$. For a hard example, we set $x_0 =
|1- y_0|$ with probability $\mu/ 2 \mu'$.  Then, we can take $\muhard
= 1/2$ and the proof of Theorem~\ref{thm:bayesinconsistent} holds
unchanged.  But now $c_0$ {\em is\/} the Bayes optimal classifier
relative to $D$, as is easy to see. Thus, Bayesian inference can be 
inconsistent for classification in {\em both\/} case (a) (no Bayes
act in $\cC$) and case (b) (heteroskedasticity).

\subsubsection{Why is the result interesting for a Bayesian?}
Here we answer several objections that a Bayesian might have to our work.
\paragraph{Bayesian inference has never been designed to work
under misspecification. So why is the result relevant?} \ 
\\
We would maintain that in {\em practice}, Bayesian inference is applied
{\em all the time\/} under misspecification in classification
problems \cite{Grunwald98b}. It is
very hard to avoid misspecification with Bayesian classification,
since the modeler often has no idea about the noise-generating
process. Even though it may be known that noise is not
homoskedastic, it may be practically impossible to incorporate 
all ways in which the noise may depend on $x$ into the prior.

\paragraph{It is already  well-known that Bayesian inference
can be inconsistent even if $\cP$ is {\em well-specified}, i.e. if it
contains $D$ \cite{DiaconisF86}. So why is our result interesting?} \ 
\\
The
(in)famous inconsistency results by Diaconis and
Freedman \cite{DiaconisF86} are based on nonparametric inference with uncountable sets
$\cP$. Their theorems require that the true $p$ has small prior
density, and in fact prior {\em mass\/} $0$ (see also \cite{Barron98}). In
contrast, Theorem~\ref{thm:bayesinconsistent} still holds if we assign
$p_{c_0,\erateD(c_0)}$ arbitrarily large prior mass $< 1$, which, by
the Blackwell-Dubins theorem, guarantees consistency if $\cP$ is
well-specified. We show that consistency may {\em still\/} fail
dramatically if $\cP$ is misspecified. This is interesting because
even under misspecification, Bayes is consistent under fairly broad
conditions \cite{BunkeM98,KleijnvdV04}, in the sense that the
posterior concentrates on a neighborhood of the distribution that
minimizes KL-divergence to the true $D$. Thus, we feel our result is relevant at least 
from the {\em inconsistency under misspecification\/} interpretation. 
\paragraph{So how can our result co-exist with theorems establishing
  Bayesian consistency under misspecification?} \ 
\\
Such results are typically proved under either one of the following
two assumptions:
\begin{enumerate}
\item The set of distributions $\cP$ is `simple', for example,
  finite-dimensional parametric. In such cases, ML estimation is
  usually also consistent - thus, for large $m$ the role of the prior
  becomes negligible. In case $\cP$ corresponds to a classification
  model $\cC$, this would obtain, for example, if $\cC$ were finite or had finite VC-dimension.
\item $\cP$ may be arbitrarily large or complex, but it is 
{\em convex\/}: any finite mixture of elements of $\cP$ is an element
of $\cP$. An example is the family of Gaussian mixtures with an
arbitrary but finite number of components \cite{Li97}.
\end{enumerate}
It is clear that our setup violates both conditions: $\cC$ has
infinite VC-dimension, and the corresponding $\cP$ is not closed under
taking mixtures. This suggests that we could make Bayes consistent
again if, instead of $\cP$, we would base inferences on its convex
closure $\overline{\cP}$. Computational difficulties aside,this
approach will not work, since we now use the crucial part (1) of
Proposition~\ref{pro:trafo} will not hold any more: the conditional
distribution in $\overline{\cP}$ closest in KL-divergence to the true
$p_{D}(y|x)$, when used for classification, may end up having larger
generalization error (expected $0/1$-loss) than the optimal classifier
$c^*$ in the set $\cC$ on which $\cP$ was based. We will give an
explicit example of this in the journal version of this paper. Thus,
with a prior on $\overline{\cP}$, the Bayesian posterior will
converge, but potentially it converges to a distribution that is
suboptimal in the performance measure we are interested in.

\paragraph{How `standard' is the conversion from classifiers to
  probability distributions on which our results are based?} \ 
\\
One may argue that our notion of `converting' classifiers into
probability distributions is not always what Bayesians do in practice.
For classifiers which produce {\em real-valued} output, such as neural
networks and support vector machines, our transformation coincides
with the logistic regression transformation, which is a standard
Bayesian tool; see for example \cite{Jordan95,Tipping01}. But our
theorems are based on classifiers with $0/1$-output. With the
exception of decision trees, such classifiers have not been addresses
frequently in the Bayesian literature. Decision trees have usually
been converted to conditional distributions differently, by assuming a
different noise rate {\em in each leaf\/} of the decision tree \cite{HeckermanCMRK00}. This makes the set of all decision trees on a given input space
$\cX$ coincide with the set of all conditional distributions on $\cX$,
and thus avoids the misspecification problem, at the cost of using a
much larger model space. 

Thus, we have to concede that here is a weak point in our analysis: we
use a transformation that has mostly been applied to real-valued
classifiers, whereas our classifiers are $0/1$-valued.  Whether our
inconsistency results can be extended in a natural way to classifiers
with real-valued output remains to be seen. The fact that the Bayesian
model corresponding to such neural networks will still typically be
misspecified suggests (but does not prove) that similar scenarios may
be constructed.

\section{Interpretation from an MDL Perspective}
\label{sec:mdl}
From an MDL Perspective, the relevance of our results needs
much less discussion: the two-part code formula (\ref{eq:cmdl}) has been
used for classification by various authors; see, e.g.,
\cite{Rissanen89,QuinlanR89} and \cite{KearnsMNR97}.
\cite{Grunwald98b} first noted that in this form, by using Stirling's
approximation, (\ref{eq:cmdl}) is essentially equivalent to MAP
classification based on the models $p_{c,\theta}$ as defined in
Section~\ref{sec:algorithms}. 
Of course, there exist more refined versions of MDL based on one-part
rather than two-part codes \cite{BarronRY98}. To apply these to
classification, one somehow has to map classifiers to probability
distributions explicitly. This was already anticipated by Meir and Merhav \cite{MeirM95} who used
the transformation described in this paper to define one-part codes. The resulting approach is closely related to the Bayesian
posterior approach $\cbayes$, suggesting that a version of our
inconsistency Theorem~\ref{thm:bayesinconsistent} still applies. Rissanen \cite{Rissanen89} considered
mapping classifiers $\cC$ to distributions $\{ p_{c,\theta^*} \}$ to a
{\em single\/} value of $\theta^*$, e.g., $\theta^* = 1/3$. As
discussed in Section~\ref{sec:technical}, a version of
Theorem~\ref{thm:bayesinconsistent} still applies to the resulting
distributions.

\paragraph{How to code hypotheses -- choice of codes and priors}
It may seem that our results are in line with the
investigation of Kearns et al. \cite{KearnsMNR97}. This, however, is
not clear -- Kearns et al. consider a scenario in which two-part code
MDL for classification shows quite bad experimental performance.  
and MDL must be `consistent'. Indeed, Kearns et al. observe that for
However, according to \cite{ViswanathanWDK99}, this is caused 
by the coding method used to
encode hypotheses. This method does not take into account the precision of
parameters involved.  In the paper \cite{ViswanathanWDK99}, a slightly
different coding scheme is proposed. With this coding scheme, MDL
apparently behaves quite well on the classification problem studied by
Kearns et al.

One may transplant the arguments of Viswanathan et al.
\cite{ViswanathanWDK99} to our setup: we can only prove inconsistency
for specific choices of the prior, corresponding to particular ways of
coding hypotheses. In practice, one usually employs
hypotheses of a different nature than we do, and one can use
properties of hypotheses such as the precision with which they are
specified to come up with `reasonable' priors/codes, which possibly do not
suffer any inconsistency problems. However, the intriguing fact
remains that {\em if\/} a probabilistic model $\cP$ is well-specified,
then under {\em very\/} broad conditions MDL is consistent
\cite{BarronC91} -- {\em under almost no conditions on the prior}. Our
work shows that if a set of classifiers $\cC$ is used (corresponding
to a misspecified probability model $\cP$), then the choice of prior
becomes of crucial importance, even with an infinite amount of data.
\paragraph{Related Work}
Yamanishi and Barron \cite{Yamanishi98,Barron91} proposed
modifications of the two-part MDL coding scheme so that it would be
applicable for inference with respect to general classes of predictors
and loss functions, including classification with $0/1$-loss as a
special case. Both Yamanishi and Barron prove the consistency (and
give rates of convergence) for their procedures.  Similarly,
McAllester's PAC-Bayesian method \cite{McAllester99} can be viewed as
a modification of Bayesian inference that is provably consistent for
classification, based on sophisticated extensions of the Occam's Razor bound, Theorem~\ref{thm:orb}.
These modifications anticipate our result, since it must
have been clear to the authors that without the modification, MDL (and discrete Bayesian MAP)
are not consistent for classification. Nevertheless, we seem to
be the first to have explicitly formalized and proved this\vspace*{-0 cm}.

\section{Acknowledgments}
The ideas in this paper were developed in part during the workshop {\em
Complexity and Inference}, held at the DIMACS center, Rutgers
University, June 2003.  We would like to thank Mark Hansen, Paul
Vit\'anyi and Bin Yu for organizing this workshop, and Dean Foster and Abraham Wyner
for stimulating conversations during the workshop.
\bibliographystyle{plain} 
\bibliography{master,MDL,peter}
\appendix
\end{document}